\documentclass[12pt,a4paper]{article}
\usepackage[T2A]{fontenc}
\usepackage[cp1251]{inputenc}
\usepackage{amsmath,amssymb,amsthm,amsfonts,amscd}

 \theoremstyle{plain}
\newtheorem{theorem}{Theorem}
\newtheorem*{main*}{Main Theorem}
\newtheorem*{theorem*}{Theorem}
\newtheorem{lemma}[theorem]{Lemma}

\newtheorem{corollary}[theorem]{Corollary}

\newtheorem{problem}[theorem]{Problem}

\theoremstyle{definition}

\newtheorem*{remark*}{Remark}
\newtheorem*{example*}{Example}

\def\Z{{\Bbb Z}}
\def\R{{\Bbb R}}
\def\RP{{\Bbb R}\!{\rm P}}

\def\A{{\bf A}}
\def\AA{\dot{{A}}}
\def\aa{\dot{{a}}}
\def\D{{\bf D}}

\def\I{{\bf I}}

\def\R{{\Bbb R}}

\begin{document}
\sloppy

\title{A generalized Kervaire Problem in stable homotopy groups of spheres}
\author{Petr M. Akhmetiev}
\maketitle

\begin{abstract}
We will give an elementary self-contained description of the Mahowald element in the stable homotopy group of spheres 
$\Pi_{2^l}$, $l \ge 3$. Using this construction, we prove that a generalized Kervaire Problem, formulated by the author in \cite{A}
is solved positively. 
\end{abstract}

\section*{Introduction}

The Kervaire Problem in stable homotopy groups of spheres 
is a problem to calculate all dimensions, for which there exist an element in the $2$-component of the stable homotopy group of sphere
 $\Pi_{n}$ with the Arf-invariant $1$. 
 By  the statement of the Hill-Hopkins-Ravenel Theorem \cite{H-H-R} this is possible only for $n=2^l-2$, $l \le 6$
 and, probably, for $l=7$.

According to an approach
\cite{E1}, \cite{A-E} 
the Kervaire Problem is reformulated as following. 
Let us assume that the integer $n$ is even and consider  the homomorphism 
$Imm^{fr}(n,1) \to Imm^{sf}(n-1,1)$,
which corresponds to an immersion 
$f: M^{n} \looparrowright \R^{n+1}$
in its regular cobordism class
the regular cobordism class 
  $g: N^{n-1} \looparrowright \R^n$
  of its self-intersection manifold.
  Then, let us consider the cobordism group of skew-framed   immersions  $Imm^{sf}(\frac{n}{2},\frac{n}{2})$ and
  the homomorphism 
\begin{eqnarray}\label{J}
J: Imm^{sf}(n-1,1) \to Imm^{sf}(\frac{n}{2},\frac{n}{2}),
\end{eqnarray}
which corresponds to an immersion
 $g$ a skew-framed immersion $(g-1,\kappa_1,\Psi)$, which is a restriction of  $g$ on the submanifold   $N_1^{\frac{n}{2}} \subset N^{n-1}$,
 this submanifold is defined as the transversal preimage of the universal submanifold 
	\begin{eqnarray}\label{1}
	\RP^{\infty - \frac{n}{2}-1} \subset \RP^{\infty}
\end{eqnarray}
by the mapping
 $\kappa: N^{n-1} \to \RP^{\infty}$,
 which represents the orienting class
 of the manifold $N^{n-1}$, assuming that $\kappa$ is transversal  along the submanifold  (\ref{1}).

Let us consider the homomorphism
\begin{eqnarray}\label{sip}
	sip: Imm^{sf}(\frac{n}{2},\frac{n}{2}) \to \Z/2,
	\end{eqnarray}
which corresponds to a skew-framed immersion a parity of its self-intersection points.
Let us consider the composition  (\ref{J}), (\ref{sip}):
\begin{eqnarray}\label{2}
\Theta=	sip \circ J: Imm^{sf}(n-1,1) \to \Z/2.
\end{eqnarray}
	
\begin{theorem}
	Kervaire Problem is solved positively for the dimension $n$, iff the homomorphism 
 (\ref{2}) is non-trivial.
\end{theorem}	
	
In  \cite{A} a generalized Kervaire Problem is formulated: 

\begin{problem}
	To calculate all dimensions $n$, for which the homomorphism
 (\ref{sip}) is non-trivial. If for a dimension $n$ the homomorphism is non-trivial, one say that the problem has a positive solution in the dimension $n$.
\end{problem}	

A generalized Kervaire Problem admits important applications. 

At the first, by
 \cite{A-C-R} a positive solution of the problem in a dimension $n$ implies that
 there exists a stably-framed immersion
 (this means that the normal bundle of the immersion is equipped with the stable framing) of a closed $n$-dimensional manifold 
 with the Browder-Eccles invariant $1$
 (the Browder-Eccles invariant is the twisted Arf-invariant, associated with the stable-framed immersion).
 This resuld is motivated by calculations of the Eccles-Wood filtration of stable homotopy groups of spheres.
 
 At the second, if we additionally assume that an immersed manifold $\zeta^{sf}_l \in Imm^{sf}(\frac{n}{2},\frac{n}{2})$, $\frac{n}{2}=2^l-1$,
  which determines a positive solution of the generalized Kervaire Problem has an integer orienting cohomology class, 
  the problem was investigated in 
  \cite{B-J-M}
 and is called the strong Kervaire Problem in dimension $\frac{n-2}{2}$. 
 From a positive solution of the strong Kervaire Problem  the Brown-Peterson Theorem
(a reference is given in \cite{B-J-M}) implies, that the Kervaire Problem in the dimension
$\frac{n-2}{2}$ has a positive solution.

The result of the paper is the following:

\begin{theorem}\label{Th1}
	A generalized Kervaire Problem
	has a positive solution for all integers
 $n=2^l-2$, $l \ge 2$.
\end{theorem}

As a conclusion, a calculation of the minimal codimension of skew-framed immersion of $n$-dimensional manifold with an odd number of self-intersection points, as well as a calculation of the minimal cohomology length of the orienting class for a skew-framed immersed $\frac{n}{2}$-manifold in the codimension $\frac{n}{2}$ are interesting.

The proof of Theorem
\ref{Th1} is based on a result by Mahowald \cite{M},\cite{E-Z}.
Because our approach is different with respect to the classical approach, even at a level of statements of the results, 
we will consider a statement of the Mahowald Theorem   (Section \ref{sec3}) and give a self-contained proof of this theorem. The main result is Theorem  \ref{Th1}, which is proved geometrically by a short construction.

\section{The element  $\eta_i^{sf}$, $i=3$,
	with the Hopf invariant $1$}\label{section1}
In the Section the first term $i=3$ of the family of immersions
 $\phi_i: N^{2^i-1} \looparrowright \R^{2^i}$,
 which represents a non-trivial element is the cobordism group $Imm^{sf}(2^i-1,1)$, $i \ge 3$ is constructed. For this element
 the Hopf invariant (the Hopf invariant $h(f)$ for an $r$-immersion
 $f: N^r \looparrowright \R^{r+k}$ $k<r$, is the number modulo 2 of  self-intersection points of a generic alteration of the following  immersion:
 $N^r \to \R^{r+k} \subset \R^{2r}$) satisfies the equations: $h(\phi)=1$. 

We follows an approach toward stable homotopy groups of spheres, which is based on 
the immersion theorem by Smale, see \cite{H}. Let us recall a definition of the cobordism groups
 $Imm^{sf}(n-k,k)$, $Imm^{st-sf}(n-k,k)$. 
 
Let $\varphi: M^{n-k} \looparrowright \R^{n}$ be an immersion, $\kappa: M^{n-k} \to \RP^{\infty}$ be a characteristic class, which is called the characteristic class of a skew-framing,  $\gamma \to \RP^{\infty}$ be the canonical line bundle over the real projective space $\RP^{\infty}$, $\Xi: \nu_{\varphi} \equiv \kappa^{\ast}(k\gamma)$ be an isomorphism of the normal bundle of the immersion  $\varphi$ with
the inverse image of $k$-dimensional bundle  $k\gamma$ over $\RP^{\infty}$, the Whitney sum of $k$ copies of the line bundle  $\gamma$. The isomorphism
$\Xi$ is called a skew-framing of the immersion  $\varphi$ in the codimension $k$.

Let us define the following Abelian group
 $Imm^{sf}(n-k,k)$ 
 as regular cobordism classes of immersion (up to regular cobordism), the group structure is determined by  disjoint union of the following triples:
  $(\varphi,\kappa,\Xi)$.
  The regular cobordism relation of triples is a standard relation, see
    \cite{A-F} for details.
    
    Let us say that an element
  $x \in Imm^{sf}(n-k,k)$ is represented by an immersion with an integer characteristic class of its skew-framing, if there exists in
  the regular cobordism class of the element a triple
  $(\varphi,\kappa,\Xi)$ for which the image of the mapping  $\kappa$ is inside $\RP^1 \subset \RP^{\infty}$.

Let $\varphi: M^{n-k} \looparrowright \R^{n}$ be an immersion,  $\kappa: M^{n-k} \to \RP^{\infty}$ be the characteristic class of its skew-framing, $\Xi^{st}: \nu_{\varphi} \oplus N\kappa^{\ast}(\gamma) \equiv \kappa^{\ast}((k+N)\gamma)$ be
a stable isomorphism of the normal bundle of the immersion $\varphi$ with the inverse image of a $k$-dimensional bundle $k\gamma$ over $\RP^{\infty}$, $N$ be a great positive integer.
Such an isomorphism $\Xi^{st}$ is called a stable skew-framing of the immersion
 $\varphi$ in the codimension $k$.

Let us define the cobordism group
  $Imm^{st-sf}(n-k,k)$
  as a set (the group structure is defined by disjoint union of triples described above)
   of regular cobordism classes of triples  $(\varphi,\kappa,\Xi^{st})$.
   The equivalent relations of triples is determined by regular cobordism classes of triples.

The Mahowald element
 $\eta^{fr}_i \in \Pi_{2^i}$ is defined using the element $\eta^{sf}_i$, constructed below, by the Khan-Priddy homomorphism :
$$ \lambda: Imm^{st}(2^i-1,1) \to \Pi_{2^i},$$
\cite{K}. 

In this section we construct the Mahowald element $\eta^{sf}_l$, $l=3$, which determines the base of the induction over the parameter $l \ge 3$. The following lemma explain a geometrical meaning of the
following in \cite{M}: the element $hh_3$ is sherical in ASS.

\begin{lemma}\label{L1}
	There exists an element
$x \in Imm^{sf}(7,1)$ with the Hopf invariant  $h(x)=1$.
\end{lemma}

\subsection*{Proof of Lemma \ref{L1}}

Obviously, the image of the generator
 $\Pi_7$ by the homomorphism   $forg: \Pi_ 7 \equiv Imm^{fr}(7,1) \to Imm^{sf}(7,1)$ into the group of skew-framed immersion
 satisfies the statement of the lemma.
 
 To our goal the base of the induction has to be proved alternatively.
 Let us construct another immersion (which satisfies special properties)
 $\phi: N^7 \looparrowright \R^8$,
which represents an element $x=[\phi]$, $h(x)=1$.

Define the immersion 
	 $\phi$ by the composition of a 2-sheeted covering $p: N^{2^i-1} \to \hat N^{2^i-1}$, $i=3$
	 and an immersion 
 $\hat{\varphi}: \hat N^{2^i-1} \looparrowright \R^{2^i}$.
 At the beginning, let us prove in this case that the Hopf invariant
  $h(\phi)$ is calculated by the formula:
\begin{eqnarray}\label{hopf}
\langle p^{2^i-1};[\hat N^{2^i-1}] \rangle = h(\phi), \quad i=3,
\end{eqnarray}
where $p \in H^1(\hat{N}^{2^i-1};\Z/2)$ is a characteristic cohomology class,
which represents the covering (we denote a covering and its cohomology characteristic class the same).

To get a proof, let us consider an immersion 
$\hat{\phi} = I \circ \hat{\varphi}: \hat N^{2^i-1} \looparrowright \R^{2^i} \subset \R^{2^{i+1}-2}$,
which to be assumed a self-transversal.
Let us consider an immersion
 $\phi'=\hat{\phi} \circ p: N^{2^i-1} \looparrowright \R^{2^{i+1}-2}$. 
 Denote by
$\phi: N^{2^i-1} \looparrowright \R^{2^{i+1}-2}$ a generic immersion, which is closed  to $\phi'$. One has to prove that  $h(\phi) \pmod{2}$ (the number of self-intersection points of $\phi$) is calculated by the formula: (\ref{hopf}). 

Obviously, self-intersection points of
the immersion
$\hat {\phi}$ give no contribution to the value $h(\phi)$, because a self-intersection point of $\hat{\phi}$ determines an even number of self-intersection points on the double covering. Therefore only self-intersection points of a deformation  $\phi \mapsto \phi'$ inside
a regular neighbourhood of the immersion
 $\hat{\phi}$ gives a contribution.
 It is clear (the reasoning is analogous to \cite{Sz}), that the number (modulo 2) of such points is calculated by the formula: 
$$\langle (w_1(\hat N) + p) p^{2^i-2};[\hat N^{2^i-1}] \rangle = h(\phi).$$ 
Because 
  $$\langle w_1(\hat N) p^{2^i-2};[\hat N^{2^i-1}] \rangle = 0,$$ 
  using the fact that one-dimensional manifold, which is Poincar\`e dual
  to the cocycle
   $p^{2^{i}-2} \in H^1(\hat{N};\Z/2)$, the expression for the Hopf invariant is given by the left hand side of the formula 
$(\ref{hopf})$.

Let us define a closed manifold 
 $N^{7} = (S^3 \times S^3) \tilde\times S^1$ as a semi-direct product   $S^3 \times S^3$ with the circle $S^1$. This manifold $N^7$ is fibred over $S^1$ with the fibre $S^3 \times S^3$. 
 Fibres $S^3 \times S^3$ are identified (by a monodromy of a generic loop on the base $S^1$ of the fibre)   by the involution
 $S^3 \times S^3 \to S^3 \times S^3$,
 which permutes the factors. 
 The oriented class of the manifold   $N^{7}$, which is denoted by $\kappa_N \in H^1(N^{7};\Z/2)$, is induced from the generator $H^1(S^1)$ by the projection $N^{7} \to S^1$.

Let us consider the following manifold:
   $\tilde{N}^{7}=(\RP^3 \times \RP^3) \tilde\times S^1$, 
   this manifold is fibred over $S^1$ with the fibre
   $\RP^3 \times \RP^3$,
   this manifold is a result of the quotient of the manifold
    $N^7$ by a pair of free involutions.  The monodromy of a fibre   $\RP^3 \times \RP^3$ along the base $S^1$
    is defined by the involution 
 $\RP^3 \times \RP^3 \to \RP^3 \times \RP^3$,
 which permutes the factors.
 Let us define an immersion
 $\tilde{\varphi}_0: \tilde{N}^{7} \looparrowright \R^{8}$, 
 the normal bundle of this immersion is a linear bundle, which coincides with the orienting bundle
 $\nu_{\tilde \varphi_0} \cong \tilde{\kappa}_N$.
 
 Let us define an auxiliary immersion 
 $\tilde \varphi: \tilde{N}^7 \looparrowright \R^9$. Consider an immersion  $f: \RP^3 \looparrowright D^4$ of the codimension $1$ and consider the Cartesian product $f_1 \times f_2: \RP^3 \times \RP^3 \looparrowright D^4_1 \times D^4_2$. Let us consider the manifold  $(D_1^4 \times D^4_2) \rtimes_T S^1$, which is defined as the semi-direct product  $D^4_1 \times D^4_2$ with the circle, along the circle the 4-disks are permuted.  
 Let us define an immersion $\tilde \varphi': \tilde{N}^7 \looparrowright (D_1^4 \times D^4_2) \rtimes_T S^1$, as
 an immersion, which is induced from the immersion $f_1 \times f_2$ by a semi-direct product. The normal bundle
 of this immersion is $2$-dimensional,
 this bundle is the Whitney sum $\tilde \kappa_N \oplus \varepsilon$.
 Let us define an immersion $\tilde \varphi$ as the composition  $I \circ \tilde \varphi'$,
 where $I: (D_1^4 \times D^4_2) \rtimes_T S^1 \subset \R^{10}$ is an embedding of a (thin) cylinder.
 
By the Hirsch Theorem \cite{H}, starting by the immersion
$\tilde{\varphi}$, one may define a compressed immersion  $\tilde{\varphi_0}: \tilde N^7 \looparrowright \R^8$,
uniquely, up to a concordance. 
Let us now define the required manifold
  $\hat{N}^7$, which is equipped by a characteristic class $p \in H^1(\hat{N}^7;\Z/2)$, which satisfies the condition: 
\begin{eqnarray}\label{hopf1}
	\langle p^{2^i-1};[\hat N^{2^i-1}] \rangle = 1 \pmod{2}, \quad i=3,
\end{eqnarray} 
and, moreover, an immersion  $\hat{\varphi}: \hat{N}^7 \looparrowright \R^8$. 

For this purpose let us consider the standard 4-sheeted covering 
$p: N^{7}=(S^3 \times S^3) \tilde{\times} S^1 \to \hat{N}^{7}$, which is induced by the pair of coverings of factors:  $S^3 \to \RP^3$.


Let us define the following manifold
 $\hat{N}^{7}$, as a manifold of pairs of two-point configurations  $[(x_1,x_2),(y_1,y_2)]$ in $N^7$, which satisfies the condition: $x_1=Tx_2$, 
$y_1 = Ty_2$, where $T$ is the transposition of points in a fibre of a corresponding covering factor:  $S^3 \to \RP^3$. The two pairs of 2-point configurations are equivalent, if one configuration is the result of the transposition of points in each 2-point configurations.
As a result, the manifold
 $\hat{N}^7$ as a base of 2-sheeted covering with the covering space $N^7$ is well-defined. Let us denote  by  $p \in H^1(\hat{N}^7;\Z/2)$ the characteristic class of this covering.

\begin{lemma}\label{L3}
	For the characteristic number, constructed by the only cohomology class 
 $p$ in $H^1(\hat{N}^7)$,  the equation $(\ref{hopf1})$ is satisfied. 
\end{lemma}

\subsection*{Proof of Lemma \ref{L3}}

Let us consider the standard $7$-sphere
 $S^7$ as the join of the two standard  3-spheres $S^3$: $S^7 = S^3 \ast S^3$.
 The antipodal involution
 $T_{S^7}$ on $S^7$ is induced from the antipodal involutions $T_{S^3}$ on the factors of the join. Consider  6-torus $P^6 \subset S^3 \ast S^3$, which is defined as the middle submanifold in the join. Obviously,  $P^6 \cong S^3 \times S^3$
 and the restriction of the involution $T_{S^7}$ on this torus $P^6$ 
 coincides with the diagonal of the two involutions $T_{S^3}$ on the factors.

Let us construct a mapping
$\hat{N}^7 \to \RP^7$ of the degree $1 \pmod{2}$, and which corresponds  with
the cohomology class $p$ by the stabilization  $\RP^7 \subset \RP^{\infty}$ of the image. Consider the projection
$\hat{N}^7 \to S^1$ on the generator and
denote by
 $\hat{K}^6 \subset \hat{N}^7$ the submanifold, which is defined as the fibre of the projection. Clearly,   $\hat{K}^6 = S^3 \times S^3/T_{diag}$, where $T_{diag}$ is the diagonal involution $T_{diag}=T_{S^3} \times T_{S^3}$. Let us consider the cutting of
 $\hat{N}^7$ along $\hat{K}^6$ and define a mapping  $\hat{F}: \hat{N}^7 \setminus \hat{K}^6 \to \RP^7$, which maps the two components of the boundary into two disjoint copies 
$\RP^3 \cup \RP^3 \subset \RP^7$ as following. Describe  $\hat{F}$ using the equivariant covering mapping $F: N^7 \setminus K^6 \to S^3 \ast S^3$, 
which is mapped the two components of the boundary onto the disjoint spheres, the generators of the join.

 By the construction, the restriction of the equivariant mapping $F$  on the first (on the second, correspondingly)
 component of the boundary coincides with the composition of the equivariant projection  
$p_i: S^3 \times S^3/ T_{diag}$ on the first $i=1$ (on the second, correspondingly)
$i=2$ coordinate with the equivariant inclusion $I_i: S^3 \subset S^3 \ast S^3$ of the first (the second, correspondingly) coordinate. 

Let us prove that the mapping
 $p_1$ is equivariant homotopic to $p_2$ after the equivariant stabilization of the images $S^3 \subset S^6$ in a common manifold. To prove this, let us define
 the total obstruction of an equivariant homotopy of the two equivariant mappings
$f_i : S^3 \times S^3 \to S^6$, $i=1,2$.
This obstruction is defined using the standard definition, this total obstruction is denoted by
  $o(f_1,f_2)$ and belongs to $\Z$.
  
  To prove this fact we recall that the involution in the preimage $S^3 \times S^3$ is defined as the diagonal antipodal involution $T_{diag}$ on the factors and this involution preserves the orientation. In the image $S^6$ the involution
  is the standard antipodal involution and inverses the orientation. An obstruction $o(f_1,f_2)$ looks like an analogous obstruction in the problem
  for the covering mappings $\bar{f}_1, \bar{f}_2: T^2 \to \RP^2$, where $T^2$ is the standard torus. This obstruction is identified with the half degree of the covering mappings in the prescribed coordinate system in the target.
  
  In the considered situation the involution
   $T_{diag}$ keeps the orientation, and the standard involution $T_{S^6}$ reveres the orientation, the integer degree $deg(f_i)$ of the two equivariant mappings are related with the obstruction by the formula:
$$ o(f_1,f_2) =\frac{1}{2} (deg(\bar{f}_1) - deg(\bar{f}_2)).$$
Because
 $deg(\bar{p}_1)=deg(\bar{p}_2)=0$, the mappings   $p_1$, $p_2$ are 
 homotopic inside
 $\RP^6$.

Let us consider the submanifold
 $\RP^6 \subset \RP^7$, which contains the images of the equivariant mappings $I_i \circ p_i$, $i=1,2$. The mappings $f_1,f_2$ are equivariant homotopic, therefore an extension of the mapping $F$
 on the manifold $N^7$ is well defined, and extra regular preimages are not added. The calculation proves Lemma 
  \ref{L3}, as well as Lemma  \ref{L1}. \qed

\section{The Mahowald elements $\eta_i^{sf}$, $i>3$ }\label{sec3}
In this section let us consider the next step of the induction $i=4$, and construct
an immersion $\phi: N^{15} \looparrowright \R^{16}$, which represents the element
 $[\phi]=\eta_4^{sf}$ in the group $Imm^{sf}(15,1)$, using the element  $\eta_3^{sf} \in Imm^{sf}(7,1)$, constructed above.

Let us repeat a reasoning and define $\tilde{M}^{15} = (\hat{N}^7 \times \hat{N}^7) \tilde \times S^1$. As the component of this semi-direct product let us take the manifold  $\hat{N}^{7}$, which is constructed in the section
\ref{section1}. Then let us define the tower of $2$-sheeted coverings:
\begin{eqnarray}\label{cov}
M^{15} \stackrel{p}{\longrightarrow} \hat{M}^{15} \stackrel{\hat{p}}{\longrightarrow} \tilde{M}^{15},
\end{eqnarray}
the bottom covering $\hat{p}$ is called
a primary covering, the top covering $p$ is called a secondary covering. To indicate the base of the covering we will use the denotation $p_{\hat{M}}$ e.t.c.  
The following  formula, which is an analogue of the formula (\ref{hopf}) in the case $i=4$ is satisfied:

\begin{eqnarray}\label{hopf2}
	\langle p^{2^i-1};[\hat M^{2^i-1}] \rangle = h(\phi), \quad i=3.
\end{eqnarray}
In the proof one can put
$\hat{N}^7 = \RP^7$, then to pass to the required case and prove that the equivariant composition 
$p \times p: N^7 \times N^7 \to S^7 \times S^7$, 
where the equivariant structure in the image and the preimage are defined by the free diagonal involutions. The mapping is equivariant homotopic inside $S^{14}$
to an equivariant mapping, which is defined by the composition of the transposition of the factors.

Then, starting by the manifold 
 $\hat{M}^{15}$ the required manifold  $\hat{N}^{15}$, equipped by the covering $p: N^{15} \to \hat{N}^{15}$, and by the immersion $\varphi=\hat{\varphi} \circ p: N^{15} \looparrowright \R^{15}$ is constructed. This immersion represents the described element  $\eta_4^{sf} \in Imm^{sf}(15,1)$ with Hopf invariant $1$. 

Let us consider the Cartesian product
 $\hat{\varphi}_7 \times \hat{\varphi}_7: \hat{N}^7 \times \hat{N}^7 \looparrowright \R^{8} \times \R^8$ of the two copies of the immersions.
 The normal bundle  $\nu(\hat{\varphi}_7 \times \hat{\varphi}_7)$ of this immersion is $2$-dimensional and is represented as the Whitney sum of the two copies of the line bundle  $\nu(\hat{\varphi})$. 

Let us consider the Cartesian product
$D^8 \times D^8$ of two open $8$-disks and the open manifold 
$\tilde{K}^{17}=(D^8 \times D^8) \tilde{\times} S^1$, where the semi-direct product is defined by the monodromy  $m: D^8 \times D^8 \to D^8 \times D^8$ along the generator $S^1$ by the formula: $m(x,y)=(y,x)$.
Evidently,  
$\tilde{K}^{17}$ is diffeomorphic to $D^8 \times D^8 \times S^1$ and there exists an embedding  $i_{\tilde{K}}: \tilde{K}^{17} \subset \R^{17}$ (the construction below does not depend on a choice of the embedding). Let us consider the manifold   
$\tilde{M}^{15}=(\hat{N}^7 \times \hat{N}^7) \tilde{\times} S^1$ (which is analogous to the manifold  $(\RP^3 \times \RP^3) \tilde \times S^1$ constructed in the section  \ref{sec1}). Let us consider an immersion $j: (\hat{N}^7 \times \hat{N}^7) \tilde{\times} S^1 \looparrowright \tilde{K}^{17}$ and the composition $\tilde{\psi} = j \circ i_{\tilde{K}}: (\hat{N}^7 \times \hat{N}^7) \tilde{\times} S^1 \looparrowright \R^{17}$. The normal bundle  $\nu_{\tilde{\psi}}$ of the immersion $\tilde{\psi}$ admits a reduction of the structured group
 $O(2)$ to the discrete subgroup $\D$. 

The dihedral group
 $\D$ contains $8$ elements and has the corepresentation $\{a,b\vert a^2 = b^4 = 1; [a,b]=b^2\}$. The element
$b^2 \in \D$ is also will be denoted denoted $-1 \in \D$. The group
$D$ is represented in  $O(2)$ as following: the element $b$ is represented by the rotation of the plane trough the angle  $\frac{\pi}{2}$, the element $a$
is represented by the symmetry with respect to the bisector of the first-third coordinate quarter.
 Let us denote $A = a \circ b$, this element is represented by the symmetry with respect to the horizontal coordinate axis. Let us denote  $\AA = b^2A=-A$, $\aa = b^2a = -a$.

Over the manifold 
$\tilde{M}^{15}$ the following tower of 2-sheeted coverings 
(\ref{cov}) is well-defined,
this tower corresponds to the following system of subgroups:
\begin{eqnarray}\label{sec1}
	 \I_a \subset \I_{a,\aa} \subset \D.
\end{eqnarray}	  
The structured group of the normal bundle
 $\nu_{\tilde{\psi}}$ of the immersion  $\tilde{M}^{15} \looparrowright \R^{17}$ is also isomorphic to the group $\D$. 
 To distinguish the two isomorphic sequence of subgroups: the sequence of coverings over the manifold $\tilde{M}^{15}$ and the sequence of reductions
 of the structured group of the normal bundle, the subgroup of coverings tower
 is denoted by 
 (\ref{sec1}), a subgroup of the structured group is denoted by the corresponding collection of elements.
 In particular, the sequence of subgroups
 (\ref{sec1}) of the covering corresponds to the following sequence of monodromy: 
\begin{eqnarray}\label{sec2}
\{a\} \subset \{a,\aa\} \subset \{a,b\} \simeq \D.
\end{eqnarray}
The element
$A$ of the residue class of the subgroup $\I_{a,\aa} \subset \D$ acts by an involution on the covering manifold $\hat{M}^{15}$ of the covering $\hat{p}$,
this action is denoted by  $T_A$. 
The element  $-1=b^2$ in the residue class of the subgroup $\I_a \subset \I_{a,\aa}$ acts on $M^{15}$ by a covering translation $M^{15} \to \hat{M}^{15}$, this action is denoted by  $T_{b^2}$. The restriction of $T_A$ on the submanifold
\begin{eqnarray}\label{P}
M^{15} \supset (N^7  \times N^7)/\{b^2\} 
\end{eqnarray}
 is defined by the formula:  $(n_1,n_2)\mapsto(T_A(n_1),n_2)$.
 
 Let us consider the subgroups
 $\{a\} \subset \{a,b\}$, $\I_{a} \subset \D$, which are not normal subgroups.  The normaliser of the subgroup $\{a\}$ is the subgroup  $\{a,\aa\} \subset \{a,b\}$. Analogously, the normaliser of $\I_{a}$ is the subgroup $\I_{a,\aa} \subset \D$, recall that $\D$ is the structured group of the bundle  $\nu_{\tilde{\psi}}$.

The subgroup $\I_{a} \subset \D$ is a subgroup of the covering
$\hat{p} \circ p: M^{15} \to \tilde{M}^{15}$,  see (\ref{cov}), which is not a  covering with a normal subgroup.  The subgroup
$\I_{A,\AA} \subset \D$  is visualized as a translation group of the covering $M^{15} \to \tilde{M}^{15}$, over the submanifold $\hat{N}^7 \times \hat{N}^7 \subset \tilde{M}^{15}$; this group determines transformations of sheets of the covering.   

Let us define the manifold
$\hat{N}^{15}$. Let us consider the vector bundle  $\nu_{\tilde{\psi}}$ over $\tilde{M}^{15}$
and denote
\begin{eqnarray}\label{tw}
\nu^{tw}_{\tilde{\psi}}
= \nu_{\tilde{\psi}} \otimes \lambda,
\end{eqnarray}
 where $\lambda = \lambda(\hat{p})$ is a linear bundle over $\tilde{M}^{15}$, which is associated with the covering $\hat{p}: \hat{M}^{15} \to \tilde{M}^{15}$. Denote by
$\tilde{L}^{13} \subset \tilde{M}^{15}$  the Euler class of the bundle
$\nu^{tw}_{\tilde{\psi}}$ (recall, the Euler class of a vector bundle over a base manifold is defined as a submanifold of zero sections of the vector bundle, in particular the normal bundle of the Euler class submanifold is isomorphic to the bundle itself, where  $\tilde{K}^{14} \subset \tilde{M}^{15}$
is the boundary of the disk-bundle associated with the vector bundle (\ref{tw}), which is the total space of the corresponding $S^1$-bundle over 
  $\tilde{L}^{13}$. For the formula (\ref{tw}) of the Euler class of the twisted section see \cite{K1}.

The total space
 $D_{\lambda}^{16}$ of the disk bundle, which is associated with the bundle 
$\lambda$, is a $16$-dimensional manifold with the boundary $\hat{M}^{15} = \partial(D_{\lambda}^{16})$. 
In the manifold
with boundary $D_{\lambda}^{16}$ let us define a subspace $\hat{N}^{15} \subset D_{\lambda}^{16}$, which is the required closed $15$-dimensional manifold.

Take the submanifold  $\tilde{L}^{13}$ 
and take its preimage by the projection $D_{\lambda}^{16} \to \tilde{M}^{15}$
on the base of the bundle. 
Let us cut-out fibres over the regular neighbourhood 
$U(\hat{L}^{14})$ in  $D_{\lambda}^{16}$.
The boundary of the obtained manifold consists of a lateral surface of the handle and the bottom of the handle, the bottom is the double covering over  $\hat{N}^{15}_{pred} \subset D_{\lambda}^{16} \to \tilde{M}^{15}$.
Denote this closed manifold by  $\hat{N}^{15}$ (this manifold is the required manifold), the decomposition of this manifold as the lateral surface and the bottom of the handle denote by

\begin{eqnarray}\label{hatN}
	\hat{N}^{15} = \hat{N}^{15}_{pred} \cup \hat{K}^{15}.
\end{eqnarray}
Let us consider the projection $q: \hat{N}^{15}_{pred} \to \tilde{M}^{15}\setminus U(\tilde{L}^{13})$, as the restriction of the projection   $D_{\lambda}^{16} \to \tilde{M}^{15}$ on the submanifold. By this projection the image  $r(\hat{K}^{15})$ coincides with $\partial(U(\tilde{L}^{13}))$, the image
 $r(\hat{N}^{15})$ coincides with the interior $\tilde{M}^{15}\setminus U(\tilde{L}^{13})$.

Let us consider the lateral surface of the handle
$\hat{K}^{15}$. 
The middle submanifold of the lateral surface $\hat{K}^{15}$ will be denoted by $\hat{K}^{14}_m \subset \hat{K}^{15}$. The natural projection 
$\hat{K}^{15} \to \hat{K}^{14}_m$ is a $D^1$-bundle; the restriction of
this projection  on the boundary $\partial(\hat{K}^{15})$ of the lateral handle is the double covering over $\hat{K}^{14}_m$. Let us denote $\partial(\hat{K}^{15})$ by $K^{14}_m$. This covering 
coincides with the pull-back $\lambda \mapsto \lambda^!$, $\lambda^! \simeq \varepsilon$, by the bottom covering in (\ref{cov}) of the spherization $S(\nu_{\hat{L}} \otimes \lambda)$ of the normal bundle of the submanifold $\hat{L}^{13} \subset \hat{M}^{15}$, which is a copy the restriction of the normal bundle $\nu_{\hat{M}}$  over its own Euler class.
By construction, there is a cohomology class:
\begin{eqnarray}\label{past}
p: K^{14}_m \to \RP^{\infty},
\end{eqnarray}
which is the restriction of the top class in (\ref{cov}).

 We may identify 
$\hat{K}_m^{14}$ with 
$S(\nu_{\tilde{L}})\otimes \lambda$, i.e. with the spherization of the normal Euler class of the bundle (\ref{tw}). Obviously, 
\begin{eqnarray}\label{Shand}
S(\nu_{\tilde{L}}\otimes \lambda)
= S(\nu_{\tilde{L}})\otimes \lambda.
\end{eqnarray}
The left hand side of the formula (\ref{Shand}) looks more naturally.
The covering $K^{14}_m$ (\ref{past}) is defined as $S(\nu_{\hat{L}} \otimes \lambda^{!})$.

On the manifold $\hat{K}_m^{14}$ the antipodal free involution is well-defined:
$$ J_{\hat{K}_m}: \hat{K}_m^{14} \to \hat{K}_m^{14},$$
the quotient $\hat{K}_m^{14}/J_{\hat{K}_m}$ is naturally identified with the projectivization  $P(\nu_{\tilde{L}}\otimes \lambda)$
corresponded with the bundle (\ref{Shand}). Let us denote
this manifold by $\tilde{K}^{14}_m$.

 Let us consider the following diagram of the 2-sheeted coverings: 
\begin{eqnarray}\label{diag}
	\begin{array}{ccccc}
 S(\nu_{\tilde{L}} \otimes \lambda) = &\hat{K}^{14}_m 	& \stackrel{\tilde{\pi}}{\longrightarrow} & \tilde{K}_m^{14} & = P(\nu_{\tilde{L}} \otimes \lambda)\\
&\hat{q}\uparrow &  &  \uparrow \tilde{q}&\\
 S(\nu_{\hat{L}} \otimes \lambda^{!}) = &K^{14}_m & \stackrel{\hat{\pi}}{\longrightarrow}  & K^{14}_m/J_{K_m} & = P(\nu_{\hat{L}} \otimes \lambda^{!}). \\
\end{array}
\end{eqnarray}

  \begin{lemma}\label{new}
  	 There exists a cohomology class
  	$\alpha: K^{14}_m \longrightarrow \RP^{\infty}$, which is 
  	 pull-back of $p$, defined by $(\ref{past})$, by the projection $K_m^{14} \to \hat{L}^{13} \subset \hat{M}^{15}$.
  	  	\end{lemma}

  	\begin{corollary}\label{new2}
  		There exist a cohomology class
  		$\hat{\beta}: \hat{K}^{14}_m \longrightarrow \RP^{\infty}$, and the pull-back by the covering $\hat{q}$ of this class  coincides with the class $\alpha: K^{14}_m \longrightarrow \RP^{\infty}$. 
  	\end{corollary}


 \subsection*{Proof of Lemma \ref{new}}
  	
  Take the following diagram, which contains the diagram 	(\ref{diag}) as the Euler class:
  \begin{eqnarray}\label{diag2}
  	\begin{array}{ccccc}
  		S(\nu_{\tilde{M}} \otimes \lambda) = &\hat{Q}^{16}_m 	& \stackrel{\tilde{\pi}}
  		{\longrightarrow} & \tilde{Q}_m^{16} & = P(\nu_{\tilde{M}} \otimes \lambda)\\
  		&\hat{q}\uparrow &  &  \uparrow \tilde{q}&\\
  		S(\nu_{\hat{M}} \otimes \lambda^{!}) = &Q^{16}_m & \stackrel{\hat{\pi}}
  		{\longrightarrow}  & Q^{16}_m/J_{Q_m} & 
  		 = P(\nu_{\hat{M}} \otimes \lambda^!)\\  		  		
  	 \end{array}
   \end{eqnarray}

  		Let us consider the class
  		$\alpha_{\hat{R}}$ on $\hat{R}^{16}_m$, which is defined by the standard characters $\I_A$ (by the projection $\I_{A \times \AA} \to \I_A)$. Let us prove the relation 
  		\begin{equation}\label{bet1}
  			T_{\hat{R}_m} \circ \alpha_{\hat{R}} = \alpha_{\hat{R}}.
  		\end{equation}

In the diagram (\ref{diag2}) spaces and coverings are defined analogously with (\ref{diag}), instead of the manifold $\tilde{L}^{13}$, the common base of the covering, we get the manifold $\tilde{M}^{15}$ equipped by the bundles described by the formula (\ref{tw}). 

Let us consider the dyagram of spaces, which are double coverings over the corresponding spaces of the diagram (\ref{diag2}):

\begin{eqnarray}\label{diag3}
	\begin{array}{ccc}
		\hat{R}^{16}_m 	& \stackrel{\tilde{\pi}}{\longrightarrow} & \tilde{R}_m^{16} \\
		\hat{q}\uparrow &  &  \uparrow \tilde{q}\\
		R^{16}_m & \stackrel{\hat{\pi}}{\longrightarrow}  & R^{16}_m/J_{R_m}, \\
	\end{array}
\end{eqnarray}
which is defined as following.

 Recall that $\tilde{M}^{15} = \hat{N}^{7} \times \hat{N}^7 \rtimes_{\tilde{T}} S^1$. The submanifold $\hat{N}^7 \times \hat{N}^7 \subset \tilde{M}^{15} $ is the ordered 2-points configuration space of $\hat{N}^7$, the involution $\tilde{T}$ permutes  points in pairs. The covering $\hat{M}^{15} \to \tilde{M}^{15}$ is described in (\ref{cov}), using the diagram (\ref{sec1}), we get $\hat{M}^{15} = (N^7 \times N^7)/T_{diag} \rtimes_{\hat{T}} S^1$,
 where, $T_{diag}$ is the diagonal of the involution  $N^7 \to \hat{N^7}$ on the Cartesian product, this involution corresponds to the element $b^2 \in \D$; the involution $\hat{T}$ permutes points in pairs. The following covering $\hat{N}^{7} \times \hat{N}^7 \times S^1 \to \tilde{M}^{15}$, which is the covering with the subgroup $\I_{A \times \AA} \subset \D$, induces the coverings of spaces of (\ref{diag3}) over corresponding spaces of (\ref{diag2}).
 
 Let us describe the cohomology classes 
 $\alpha_{R}: R_m^{16} \to \RP^{\infty}$. Then let us prove that the following formula is satisfied: 
 \begin{eqnarray}\label{bet}
 \alpha_{R} = T_R \circ \alpha_{R},
 \end{eqnarray}
  where $T_R$ is the involution of the covering $R^{16}_m \to Q^{16}_m$.
 
 Define the class $\alpha_{R}$ as the pull-back of the class  $[T_{diag}]^{\ast} \in H^1((N^7 \times N^7)/T_{diag})$  by the projection $R_m^{16} \to (N^{7}\times N^{7})/T_{diag} \times S^1 \to (N^{7}\times N^{7})/T_{diag}$. 
 The fact (\ref{bet}) is a corollary that $\alpha_R$ coinsids to the standard $\I_{ \aa}$-character (is invariant with respect to $\I_a$-transformation).
 The cohomology class $\alpha_{Q} \in H^1(Q_m^{16})$ is well-defined.

 Lemma \ref{new} is proved. \qed

   \subsection*{Proof of Lemma \ref{new2}}
   
  Take the left-bottom manifold $K_m^{14}$ in (\ref{diag}) and
  consider the section $\varepsilon = \lambda^!$  as the pull-back of $\lambda$ of the bundle
  $\nu_{\tilde{L}} \otimes \lambda$,
  restricted
  outside the Euler class $\tilde{L}^{13} \subset \tilde{M}^{15}$ of the  cross-section. Below we will speak of
  the section  
  $\lambda^!$ as a section of the bundle $\nu_{\hat{L}}$ instead of a section
   of the bundle  $(\nu_{\hat{L}} \otimes \lambda)^!$, which is equivalent.

   An open neighbourhood of $K^{14}_m \subset Q^{16}_m$ is referred as the regular domain in $Q^{16}_m$.
   A regular domain of $K^{14}_m \subset Q^{16}_m$ and its preimage in $R_m^{16}$  are denoted by
 $Q^{16}_{m;reg}$, $R^{16}_{m;reg}$ correspondingly.  
  
 Let us denote by $H^{+}(X) \subset H_1(X)$ for an non-orientable closed manifold $X$ (we take below $X=K^{14}_m$)
 the index 2 subgroup, which is the kernel of the orientation homomorphism:
 $Ker (w): H_1(X) \to \Z/2$. Analogously, define the non-trivial resedue class of the orientation homomorphism by $H^-(X)$.
 Obviously, $H^1(X) = H^+(X) \cup H^-(X)$.
    Denote by 
    \begin{eqnarray}\label{oplus}
 H^{\oplus}(K^{14}_m) = [H^+(K^{14}_m) \cap Ker(\I_{a})]
 \cup [H^-(K^{14}_m) \cap \overline{Ker}(\I_{a})],
 \end{eqnarray}
 where $\I_{a}$ is the homomorphism
 $H_1(K^{14}_m) \to \I_{a \times \aa} \to \I_{a}$. 
 This denotation $ H^{\oplus}(K^{14}_m)$ means that we take in $H_1(K^{14}_m)$ the kernel of the homomorphism $A+\AA: \{a,b\} \to \Z/2$, described using the normal bundle structure. By this convention, one may use analogous denotations for different manifolds in the left collomn in the diagrame (\ref{diag}), and in the diagrames (\ref{diag2}), (\ref{diag3}) and for   
 submanifolds $R_{m;reg}$, $Q_{m;reg}$, defined below. 

  For an arbitrary homology class $[l]_{\ast} \in H^{\oplus}(K^{14}_m)$ let us
  consider a connected loop, which represents the same class on $R^{16}_{m;reg}$, and  the double covering $l_R$ over the loop in $R^{16}_{m;reg}$ (the left-bottom space in (\ref{diag3})).
 The loop $l_R$ is visualized as  a two-component closed loop (Case I, the srandard $a$-character (recall that the $a$-character coinsids to the standard character $\I_a$ ) on $l$ is trivial) and the covering $l_R$ over $l$ is disconnected, or, as a single connected loop (Case II, the $\I_a$-character  on $l$ is non-trivial). In  Case II, we may see
 $l_R$ as a segment on $R^{16}_{reg;m}$ with
 involutive end-points by $T_{R_m}$. Moreover, in the case the class $l_{\ast}$ is represented by a segment $l_R$, such that
  end-points of the segment are $T_{R_m}$-involutive, or, the projection of $l_R$
  on the circle fibre are $\I_a$-involutive. We assume that  end-points of $l_R$ are not on the section.
 In Сase I  components of the loop $l_R$ are ordered, and the numbers defined below do not depend on an order; In Case II segments are also ordered. 

     Denote by $\sharp (l) \in \{0,1\}$ the number of intersection points of
   the first component of $l_R$ with the positive base vector of cross-section $\lambda^!$. Take the involution $T_{\hat{q}}$ in the covering by the left vertical arrow of (\ref{diag}) and look the lift of $T_{\hat{q}}$ as an involution
   	$T_{\hat{q};R_{m;reg}}$ on $R^{16}_{m;reg}$ inside the regular domain.
   Let us prove that numbers $\sharp (l)$, $\sharp(T_{\hat{q};R_{m;reg}}(l))$ coincide in Case I, assuming the  the standard $A+ \AA$-character on $l$ is trivial, or, equivalently, $[l] \in H^{\oplus}(K^{14}_m)$, see the formula (\ref{oplus}).
  In the Case $II$ the numbers $\sharp (l_R)$, $\sharp(T_{\hat{q};R_{m;reg}}(l_R))$  
   are different and $[l]_{\ast} \in H^-(Q^{16}_{m;reg})$. In the both cases we calculate the indexes for components (segments) with a common number.

   Let us consider the Case I, recall
    $l \in H^+(R^{16}_{m;reg})$. 
   The bundles $\kappa_A$, $\kappa_{\AA}$ over $l_R$ are trivial (Subcase 1), or, non-trivial (Subcase 2) simultaneously. Take a coordinate system at the initial point ot the path. The section $\varepsilon$ over each component of $l_R$ is  
   projected on the spherical fibre, which is the projection on the circle
   in the subcase 1.
  The number $\sharp(l_R)$ for the first component depends on a number of rotation of the section, but the number  $\sharp(T_{\hat{q};R_{m;reg}}(l_R))$ is also depended, because the sections are antipodal.

  In the subcase 1 of the trivial bundles the section $\lambda^!$ is represented
  by a constant section projected on the circle fiber to the point  $e_1$. The involution $T_{\hat{q};R_{m;reg}}$ transforms the section into its antipodal point $-e_1$ on the circle. 
    In the subcase 2,
    we will see a closed loop as a twisted cylinder, where "twisted" means a rotation
    of its generator trough the angle $\pm \pi$ with respect to the coordinates in $S(\kappa \oplus \kappa_{\aa})$. In this case the section is translated by 
    $T_{\hat{q};R_{m;reg}}$ into the antipodal twisted cylinder.

   Assume a component $l_R$ is represented by a curve
   at the point $e_1$ (the base vector of the bundle $S(\nu_{\hat{L}} \otimes \lambda^!)$).
   Take the involution $T_{\hat{q};R_{m;reg}}$ and take the images of
   $l_R$ and of the section $\lambda^!$ by the involution. The curve  $T_{\hat{q};R_{m;reg}}(l_R)$ is a constant at the point $e_1$,  the section $T_{\hat{q};R_{m;reg}}(\lambda^!)$  is a constant at the point $-e_1$,  (recall, the bundle $\nu_{\hat{M}} \otimes \lambda^!$ is twisted by the involution $T_{\hat{q}}$).
    For the  curves  in the example we get $\sharp(l_R)=0$;
   $\sharp(T_{\hat{q};R_{m;reg}}(l_R))=0$. For an arbitrary $l_R$ the case 
   $\sharp(l_R)=1$;
  $\sharp(T_{\hat{q};R_{m;reg}}(l_R))=1$ 
   is also possible.
   
   As the conclusion, in the Case $I$, subcase 1, we get:
   $$\sharp(l_R)=\sharp(T_{\hat{q};R_{m;reg}}(l_R)).$$ 
   The case $I$, subcase 2 is analogous.
   
   Let us consider the Case $II$. We may
   assume that the bouth standard $A$, $\AA$-characters on $l_R$ are trivial
   (this case is analogous to the subcase 1 in Case $I$; the case corresponds to the condition $[l]_{\ast} \in H^{\oplus}(K^{14}_m)$, see the formula (\ref{oplus}); the case when the both characters are non-trivial is analogous  to the subcase 2, Case I). The curve $l_R$ is represented by a constant segment at the point $e_1+e_2$, the section $\lambda^!$ is represented by a constant at the point $e_1+e_2$. 
   The curve  $T_{\hat{q};R_{m;reg}}(l_R)$ is a constant at the point $e_1+e_2$,  the section $T_{\hat{q};R_{m;reg}}(\lambda^!)$  is a constant at the point $-e_1-e_2$.
  A generic alteration  gets:
   $\sharp(l_R)=1$;
   $\sharp(T_{\hat{q};R_{m;reg}}(l_R))=0$.
 As the conclusion, in the Case $II$ we get:
 $$\sharp(l_R)=\sharp(T_{\hat{q};R_{m;reg}}(l_R))+1.$$

   The calculation above means that $\sharp: H^+(Q_{m;reg}^{16}) \to \Z/2$ is a homomorphism, which is
   an $\I_a$-character. This exactly means that on an arbitrary class in $H^{\oplus}$ (see the denotation below) the  homomorphism $\sharp$ changes its value  with respect to the
   involution $T_{\hat{q}}: H_1(Q_{m;reg}^{16}) \to H_1(Q^{16}_{m;reg})$, iff the class  is not in the kernel of the homomorphism $\I_{a}$.

Obviously, the standard character $\I_{a}$  on $H^{\oplus}(K^{14}_m)$,  
admits
the same property. This implies that the homomorphism
$\sharp + \I_a: H^{\oplus}(R_{m;reg}^{16}) \to \Z/2$ is  invariant with respect $T_{\hat{q}}: H^{\oplus}(R_{m;reg}^{16}) \to H^{\oplus}(R_{m;reg}^{16})$, this homomorphism is a cohomology class in  $H^{\oplus}(\hat{R}_{m;reg}^{16})= H^{\oplus}(\hat{K}_m^{14})$,
where by $H^{\oplus}(\hat{K}_m^{14})$
the kernel of the homomorphism 
$A+\AA: H_1(\hat{K}_m^{14}) \to \Z/2$ is defined.

As the result, there exists a homomorphism
$\hat{\beta}: H_1(\hat{K}^{14}_m) \to \Z/2$, which is defined as the canonical extension of $\sharp + \I_a$.
An extension is defined by two possible ways, let us take $\hat{\beta}$ with a regular condition, using $\alpha$, described below.

Take the pull-back $\alpha = \hat{q}^{\ast}(\hat{\beta})$ of $\hat{\beta}$ by the covering $\hat{q}$.
The restriction of the cocycle $\sharp$ on $H^{\oplus}(K^{14}_m)$ is extended to the cocycle on $H^{\oplus}(\hat{M}^{15} \setminus U(\hat{L^{13}}))$, because the section $\lambda^{!}$ on this manifold
	$\hat{M}^{15} \setminus U(\hat{L^{13}})$ is well-defined. To calculate the pull-back $\sharp^!$ of $\sharp$ by the covering $\hat{q}$ take a generic intersection of the section on
		$\hat{M}^{15} \setminus U(\hat{L^{13}})$ with itself. Obviously, the extension $\sharp^!$ coincides with $w_1(\hat{M}^{15} \setminus U(\hat{L^{13}}))$, which is trivial on $H^+(\hat{M}^{15} \setminus U(\hat{L^{13}}))$ and is non-trivial on  $H^-(\hat{M}^{15} \setminus U(\hat{L^{13}}))$. The section $\lambda^!$, restricted on $\hat{M}^{15} \setminus U(\hat{L^{13}})$, is generic with respect to itself, the number of the intersection of the section with a generic alteration over a homology class $[l]_{\ast} \in H_1(\hat{M}^{15} \setminus U(\hat{L^{13}}))$ represents the number $\sharp^!(l)$. One gets the case $\sharp^!(l)=1$ when
		the bundle $\nu_{\hat{M}}$ over $l$ is non-trivial. This exactly means that $\hat{M}$ over $l$ is non-oriented. Because the assumption
		$[l]_{\ast} \in H^{\oplus}(K^{14}_m)$, or, equivalently, $[l]_{\ast} \in H^{\oplus}(\hat{M}^{15} \setminus U(\hat{L^{13}}))$, the condition $\sharp^!(l)=1$ iff $\I_a(l)=1$, or, equivalently, iff $[l]_{\ast} \in H^{\oplus}(\hat{M}^{15} \setminus U(\hat{L^{13}}))$.
		This means that on $H^{\oplus}(\hat{M}^{15} \setminus U(\hat{L^{13}}))$
		we get: $\sharp^! = \I_{a}+\I_{\aa}$.

An	extension of a cocycle from  $H^{\oplus}( \hat{K}^{14}_m)$ to $H_1(\hat{K}^{14}_m)$, say, the cocycle $\alpha$, determines the only corresponding extension of $\hat{\beta}$ (a value  $\hat{\beta}(t)$, $t \in H_1(\hat{K}^{14}_m) \setminus H^{\oplus}( \hat{K}^{14}_m)$ determines the homomorphism $\hat{\beta}$ on $H_1(\hat{K}^{14}_m)$, which coincides with the prescribed lift $\hat{\beta}$ 
on the codimension 2 subgroup $H^{\oplus}(\hat{K}_m^{14})$, where the character $A+\AA$ is trivial). Because $\alpha = \sharp^! + \I_{a} = \I_{\aa} = p$ on 	$H^{\oplus}(\hat{M}^{15} \setminus U(\hat{L}^{13}))$, the cocycle $\hat{\beta}$ is well-defined with the condition $\hat{q}^{\ast}(\hat{\beta})=\alpha$.	     Lemma \ref{new2} is proved. \qed

 To get a construction of $\eta^{sf}_4$ one may start with  a simplification, that
 the Euler class $\tilde{L}^{13}$ of $\nu_{\tilde{M}}$ is empty.
 With this assumption we may see that $\hat{M}^{15}=\hat{N}^{15}$ is a boundary of the total space of the disk bundle associated with $\lambda$ (looks like a generalized M$\ddot{\rm{o}}$bious band),
 immersed into $\R^{16}$, and therefore the immersion $g: \hat{M}^{15} \looparrowright \R^{17}$ admits a section by vector inside the band). This gives the required element $\eta^{sf}_4 \in Imm^{sf}(15,1)$ by the composition $M^{15} \to \hat{M}^{15} \looparrowright \R^{17}$, analogously with Section \ref{section1}.

Let us consider a general case. The manifold
 $\hat{N}^{15}$ (\ref{hatN}), which is equipped with the characteristic class  $p_{\hat{N}} \in H^1(\hat{N}^{15})$ by Lemma \ref{new2} is well-defined. Let us prove that  $p^{15}_{\hat{N}}$ is the fundamental class of this closed manifold. By the construction, 
$\hat{N}^{15}$ (see the formula (\ref{hatN})) is the result of a surgery of  $\hat{M}^{15}$ by a generalized handle with a $14$-dimensional bottom.
The bottom of the lateral surface of the handle is equipped with the class   $p_{\hat{K}}$, which is the pull-back of a class $\alpha$ by the projection of the lateral surface onto the middle submanifold  $\hat{K}^{14}_m$, as it is proved in Lemma \ref{new2}.
  Therefore the support of the cohomology class   $p_{\hat{K}}$ on the lateral handle $\hat{K}^{15}$ is induced by a mapping into  $14$-skeleton $\RP^{14} \subset \RP^{\infty}$.  
  As a result, the degree of the characteristic mapping
  into
   $15$-skeleton for the class of the mapping $p_{\hat{M}} \mapsto p_{\hat{N}}$ by the surgery $\hat{M}^{15} \mapsto \hat{N}^{15}$ remains fixed,   and the  characteristic number of $p_{\hat{N}}$ and of $p_{\hat{M}}$ coincide. The required equation 
   \begin{eqnarray}\label{required}
    \langle p_{\hat{N}}^{15};[\hat{N}^{15}] \rangle =1
    \end{eqnarray}
     is deduced from (\ref{hopf2}).   	
 	 
An immersion
 $\phi_{\hat{N}}: \hat{N}^{15} \looparrowright \R^{17}$ is well-defined, using $p_{\hat{N}}$ with the equation (\ref{required}) and is determined
 the required element in $Imm^{sf}(15,1)$, which will be denoted by  $\eta^{sf}_{4}$. Let us repeat reasoning from the section
  \ref{section1} and prove that the element
  $\eta^{sf}_{4}$ has the non-trivial Hopf invariant.

The construction of the elements
$\eta^{sf}_l \in Imm^{sf}(2^l-1,1)$, $l \ge 5$ is completely analogous to the case $l=4$. 

The Mahowald element $\eta^{fr}_4 \in Imm^{fr}(16,1)$ is well-defined.
It is not difficult to prove, using the
calculation of the secondary cohomology operation by the Khan-Priddy transfer
(see \cite{E1}), that the element  
$\eta^{fr}_4$ is detected by the corresponding secondary cohomology operation, described in \cite{M}.

\section{The Hopf invariant one element $\zeta^{sf}_l \in Imm^{sf}(2^l-1,2^l-1)$}

Let us prove Theorem \ref{Th1}.
Consider an immersion 
$\hat{\phi}: \hat{N}^{2^l-1} \looparrowright \R^{2^l}$, where $p \in H^1(\hat{N}^{2^l-1})$, $\langle p^{2^l-1};[\hat{N}]\rangle = 1$, $l \ge 3$, such an immersion is constructed in the section \ref{sec3}. Using this data we construct  by a regular surgery a self-transversal immersion  $\phi_2: N_2^{2^l-1} \looparrowright \R^{2^{l+1}-2}$
with an odd number of self-intersection points and a stable skew-framing  $\Xi_2: \nu_{\varphi_2} \oplus  n\kappa \simeq (2^l-1+n)\kappa$, where $\kappa$ 
is the line orienting bundle over
$N_2^{2^l-1}$. Then we calculate 
the obstruction to reduce the stable skew-framing $\Xi_2$ to a skew-framing
and we shall prove that there exists 
a modification of $(\varphi_2,\Xi_2) \mapsto (\varphi_3,\Xi_3)$
of stably skew-framed immersions, 
such that the stably skew-framed immersion $\phi_3$ is a skew-framed immersion. The surgery preserves the set of self-intersection points of $\phi_2$. This proves the theorem.

Let us consider the characteristic mappings
$p: \hat{N}^{2^l-1} \to \RP^{\infty}$
$\hat{\kappa}: \hat{N}^{2^l-1} \to \RP^{\infty}$, where
$\hat{\kappa}$ is the orienting cohomology class of the manifold  $\hat{N}^{2^l-1}$, $p$ is a classifying mapping of the given cohomology class. Let us denote by
\begin{eqnarray}\label{mu}
\hat{\mu}: \hat{N}^{2^l-1} \to \RP^{\infty}
\end{eqnarray}
a classifying mapping with the characteristic class  $\hat{\mu}=\hat{\kappa} + p$.
The mapping $\hat{\mu}$ is the classifying mapping of the line bundle
 $\hat{\kappa}\otimes p$ over $\hat{N}^{2^l-1}$. Consider the submanifold  $\hat{K}^{2^l-3} \subset \hat{N}^{2^l-1}$, which is defined as the regular preimage by $\hat{\mu}$ of the standard codimension $2$ submanifold: $\RP^{\infty-2} \subset \RP^{\infty}$. The manifold $\hat{K}^{2^l-3}$ is assumed to be connected. Denote by $U(\hat{K}^{2^l-3})$ a thin regular neighbourhood of this manifold in $\hat{N}^{2^l-1}$.

A generic immersion by a small generic (which is special near $\hat{K}^{2^l-3}$) alteration of the mapping
\begin{eqnarray}\label{phi1}
\phi_1 = I \circ \hat{\phi} \circ p: N^{2^l-1} \looparrowright \R^{2^{l+1}-2},
\end{eqnarray}
 where $p: N^{2^l-1} \to \hat{N}^{2^l-1}$ is the covering with the prescribed cohomology class, 
$\hat{\phi}: \hat{N}^{2^l-1} \looparrowright \R^{2^l}$
 is given, $I: \R^{2^l} \subset \R^{2^{l+1}-2}$ is the standard embedding. Define a regular surgery
 $\phi_1 \mapsto \phi_2$ using a generalized handle, which is defined as the disk-line bundle over the manifold $U(\hat{K}^{2^l-3})$, associated with the covering $p_K: K^{2^l-3} \to \hat{K}^{2^l-3}$. 
 The handle is taken along with the generator along a line bundle with the class $p_K$.
 As a result, we shall prove, that the immersion  
  $\phi_2$ become stably skew-framed in the ambient space $\R^{2^{l+1}-2 + C_{2^l-1}-(2^{{l+1}-2}-(2^{l}-1)+1)} = \R^{(2^l-1) + C_{2^l-1}+1} $, where $C_{2^l-1}$ is the Hurwitz number (the minimal power of $2$, such that the bundle $C_{2^l-1} \gamma \to \RP^{2^l-1}$ is isomorphic to the trivial bundle,  where $\gamma$ is the canonical line bundle over the standard skeleton $\RP^{2^l-1}$). Let us describe the surgery $N^{2^l-1} \mapsto N_2^{2^l-1}$.

Consider  the immersions
 $\hat{\phi}_1=I \circ \hat{\phi}: 
\hat{N}^{2^l-1} \looparrowright \R^{2^{l+1}-2}$,
 $\hat{\phi}^{st}_1=I_1 \circ I \circ \phi: 
\hat{N}^{2^l-1} \looparrowright \R^{2^{l+1}-2} \subset \R^{2^{l}+ C_{2^l-1}}$. Over the immersed manifold  $\hat{\phi}_1$ let us consider the trivial bundle  $C_{2^l-1}\varepsilon$, which is the restriction of the normal bundle of the embedding  $\R^{2^l}\subset \R^{2^l+C_{2^l-1}}$ 
and let us consider the direct first factor of the normal bundle $C_{2^l-1}p \subset \nu_{\phi_1^{st}}$ by the isomorphism
$C_{2^l-1}\varepsilon \simeq C_{2^l-1}p$, restricted over the submanifold $\hat{K}^{2^l-3} \subset \hat{N}^{2^l-1}$.
Without loss of a generality, let us assume, that the first factor line  bundle $p$ of the normal bundle of $\hat{\phi}_1^{st}$  is inside the hyperspace $\R^{2^{l+1}-2} \subset \R^{2^l+C_{2^l-1}}$.  

Let us consider a deformation of the immersion
 $\phi_1$ into a self-transversal immersion (without a change of the denotation) 
  (as it is proved above with the odd number of self-intersection points, each self-intersection is outside $\phi_1(U(K^{2^l-3}))$),
  which over   $\hat{K}^{2^l-3}$ coincides with a small deformation along the fibre of the line bundle $p$.
  The result of the surgery  by the 
  handle is the transformation 
  \begin{eqnarray}\label{surgery}
  \phi_1(N^{2^l-1}) \mapsto \phi_2(N^{2^l-1}_2).
  \end{eqnarray}
   We glue the disk bundle, associated with   $p$, over the boundary of the disk normal bundle of the codimension 2 submanifold $i_{\hat{K}}: \hat{K}^{2^l-3} \subset \hat{N}^{2^l-1}$. The stabilization of the immersion $\phi_2$ in the space 
  $\R^{2^l+C_{2^l-1}}$  let us denote by 
  $\phi_2^{st}$.

  The handlebody  consists the boundary, which consists of  the lateral handle, denoted by $Q^{2^{l-1}}$ (the lateral part of the boundary of the body) and the bottom. The lateral  handle contains the boundary, which is a manifold of the dimension $2^l-2$, this manifold is connected because $\hat{K}^{2^l-3}$ is connected and the restriction $p\vert \hat{K}^{2^l-3}$ is non-trivial. The  boundary $\partial(Q^{2^l-1})$ of the lateral handle is identified with the boundary of the regular neighbourhood of the submanifold   $i_{K}$, the double covering over $i_{\hat{K}}$.
  
  The lateral handle has the normal bundle
    is equipped with the line bundles $p$, $\kappa$, because the lateral handle is projected onto $U(\hat{K}^{2^l-3})$, the pull-back of  the line bundle $\hat{\kappa}$ over the base $\partial(U(\hat{K}^{2^{l}-3}))$, is denoted by 
		$\kappa$. Moreover, the terms $p$ and $\hat{\kappa}$ over the lateral  handle are isomorphic: this fact is satisfied, because $\partial(U(\hat{K}^{2^l-3}))$ admits a foliation across the meridians and the both classes are defined as rotation classes over this foliation trough the angle $\pi \pmod{2\pi}$. The  normal bundle of the lateral handle in $\R^{2^l+C_{2^l-1}}$ is the following:  
    \begin{eqnarray}\label{formula}
    (C_{2^l-1}-1)p \oplus \varepsilon \oplus \kappa ,
    \end{eqnarray}
where the term $\varepsilon$ is the tangent vector to the handlebody, which is perpendicular to the lateral handle.
       The  factor $\varepsilon$ in the formula (\ref{formula}) can be replaced into $p$ if we assume that the boundary conditions over end-points of segments, the generator of the handle, is changed into the opposite.  The formula (\ref{formula}) is modified as the following:  
     \begin{eqnarray}\label{formula3}
     	(C_{2^l-1}-1)p \oplus p \oplus \kappa.
     \end{eqnarray}
 The formula (\ref{formula3}) looks more natural with respect to the formula (\ref{formula}), because the boundary condition
 on the bottom of the lateral handle for (\ref{formula3}) is satisfied. The formula (\ref{formula3}) describes the normal bundle of the immersion $\phi_2^{st}$.

Let us modify the formula (\ref{formula3}).
Recall that the bundle $\kappa$ is the pull-back of the line bundle $\hat{\mu}$ over $\hat{N}^{2^{l}-1} \setminus U(\hat{K}^{2^l-3})$, described by (\ref{mu}). The target bundle $\hat{\mu}$
admits an integer characteristic classes, in particular, over $\hat{N}^{2^{l}-1} \setminus U(\hat{K}^{2^l-3})$ the following equation
\begin{eqnarray}\label{mueps}
 2\hat{\mu} = 2 \varepsilon
 \end{eqnarray}
  is satisfied. 

Let us do the following change:
\begin{eqnarray}\label{change}
	C_{2^l-1}p \oplus \kappa \mapsto (C_{2^l-1}+1)\kappa,
\end{eqnarray}
because $C_{2^l-1}p$, $C_{2^l-1}\kappa$ are isomorphic as the trivial bundles.
The  skew-framing of the immersion $\phi_2^{st}$ is given by the formula:
 $\Psi_2: \nu_{\phi_2^{st}} = 
(C_{2^l-1}+1)\kappa$. 
Let us prove that the skew-framing 
$\Psi_2$ is extended into a stable skew-framing of the immersion $\phi_2: N_2^{2^l-1} \looparrowright \R^{2^{l+1}-2}$ inside  $\R^{2^l+C_{2^l-1}}$. 

Using the formula (\ref{mueps}) and the modification (\ref{change}), we may rewrite the formula (\ref{formula3}) into the following: 
\begin{eqnarray}\label{eq}
\Psi_2: \nu_{\phi_1} \oplus (C_{2^l-1}+1-(2^l-1))\varepsilon = (C_{2^l-1}+1)\kappa.
\end{eqnarray}
 Because the bundle $\kappa$
 is the pull-back of  bundle $\hat{\mu}$ over
$\hat{N}^{2^l-1} \setminus U_{\hat{K}}$, analogously to (\ref{mueps}), one gets:  $(C_{2^l-1}-2^l+2)\varepsilon 
\simeq (C_{2^l-1}-2^l+2)\kappa$.

After we apply this isomorphism to the right-hand side of the equation
  (\ref{eq}), we get the following formula for the stable skew-framing $\Psi_2^{st}$ of the of the codimension $2^l-1$ immersion $\phi_2$:
\begin{eqnarray}\label{eq1}
\Psi_2^{st}: \nu_{\phi_2} \oplus (C_{2^l-1}+1-(2^l-1))\kappa = (C_{2^l-1}+1)\kappa.
\end{eqnarray}.

Let us calculate the total obstruction of
destabilization of the stable skew-framing
$\Psi_2^{st}$ into a skew-framing
 $\Psi_2$. The total obstruction is an integer, denoted by  $o(\Psi_2)$, this obstruction is, in particular, investigated in  \cite{A-C-R}. Let us recall the definition of the total obstruction specified for the connected stable skew-framed manifold $N_2^{2^l-1}$.  
 
 Let us consider an isotopy of the family of skew-sections
 $\Psi_2^{st}$ to a special   $C_{2^l-1}+1-(2^l-1)$-family
 of skew-sections. In such a special family all sections of the familly in the codimension $2^l-2$, except, probably, the only skew-section in the codimension $2^l-1$ in the left hand side of the equation
(\ref{eq1}) coincides with  corresponding skew-sections (which are considered as skew-sections of a standard skew-framing) in the right-hand side of the equation. 
An integer obstruction 
 $o(\Psi_2)$ to deform  the last considered skew-section in the right-hand side of the equation into the corresponding skew-section in the left-hand side of the equation by a deformation, which is fixed on subfamilly of special skew-sections is well-defined as following.
 
 Take the orienting covering over $N^{2^l-1}_2$, denoted by $\bar{N}^{2^l-1} \to N^{2^l-1}_2$ with the class $\kappa$ (which is a characteristic class of the skew-framing $\Psi^{st}_2$). By the assumption 
 an equivariant mapping $\bar{\Psi}: \bar{N}^{2^l-1} \to S^{2^l-1}$ is well-defined, this mapping is given by the last vector of a special skew-section. The degree of this equivariant mapping defines the obstruction $o(\Psi_2)$. The obstruction is well-defined when the orientation on $\bar{N}^{2^l-1}$ is fixed. The change of an orientation transforms the obstruction $o(\Psi_2)$ into the opposite.

Let us prove that
\begin{eqnarray}\label{o}
o(\Psi_2) = 0\pmod{4}. 
\end{eqnarray}
Denote  by $\psi_2$ a skew-section, which is not special, for example, the last skew-section of $\Psi_2$. Let us prove that the algebraic number of transversal intersection points of this skew-section $\psi_2$
with the subbundle (the marked subbundle) of the codimension $2^l-1$, generated by skew-sections in the right-hand side of the formula satisfies the equation (\ref{o}).

Take the manifold $\hat{N}$ and the immersion 
$$\hat{\phi}: \hat{N}^{2^l-1} \looparrowright \R^{2^l} \subset \R^{2^{l+1}-2} \subset \R^{(2^l-1) + C_{2^l-1}+1} $$
 by (\ref{phi1}) with its stabilization in the target.
Take a  skew-section $\hat{\psi}_2$ over $\hat{N}^{2^l-1}_2$ with the class $\hat{\mu}$, restricted on $\hat{N}^{2^l-1} \setminus \hat{K}^{2^l-3}$; consider this skew-section as a skew-section of the subbundle defined by the orthohonal complement of the inclusion $\R^{2^l} \subset \R^{2^{l+1}-2}$.
	The skew-section $\hat{\psi}_2$ is not non-degenerated, the kernels of the skew-section is a curve on $\hat{N}^{2^l-1}$, which is dual to $\hat{\mu}^{2^{l}-2}(\hat{N}^{2^l-1}) = 0$, the cohomology orienting class	$w_1(\hat{N}^{2^l-1})$ of this curve is trivial. This implies that a small horizontal generic alteration the section $\hat{\psi}_2$ does not intersect the   complement of the inclusion $\R^{2^{l+1}-2} \subset  \R^{(2^l-1) + C_{2^l-1}+1}$. 


 The speсial skew-framing $\Psi_2$  is a pull-back of a corresponding
special skew-framing $\hat{\Psi}_2$ over $\hat{N}^{2^l-1} \setminus U(\hat{K}^{2^l-3})$ (the only skew-section, along which the surgery is defined, is complicated,  can be omitted, because this skew-section is orthogonal to the tangent space of $\hat{\phi}_2$ and to the last vectors of the skew-framing $\hat{\Psi}_2$).  The skew-framing $\hat{\Psi}_2$ is defined by the standard reparametrization of the bundle
$C_{2^l-1}p$ into $C_{2^l-1}\hat{\kappa}$,  the formula (\ref{change}) is the pull-back by the construction over the base of the double covering. The equation (\ref{mueps}) is satisfied over $\hat{N}^{2^l-1} \setminus U(\hat{K}^{2^l-3})$ and the orthogonal complement  of the inclusion $\R^{2^{l+1}-2} \subset \R^{2^l+C_{2^l-1}}$ is represented by a skew-framing, the pull-back of this skew-framing over $N^{2^l-1}_2$ is used to defined the obstruction in the formula (\ref{o}).  By  the standard deformation of the skew-framing $\hat{\Psi}_2$ into a special skew-framing with prescribed skew-sections in the codimension $2^l$,  a parity of intersection points of the next skew-section
 with the marked subbundle, given by the complement of $\R^{2^{l+1}-2} \subset \R^{2^l+C_{2^l-1}}$, is predicted by a position of this section at the initial state  $\hat{\Psi}_2$ of the deformation.

For an arbitrary special skew-framing the parity of points in this intersection set is calculated using $\hat{\psi}_2$, because an arbitrary two skew-sections of the stabilized normal bundle are isotopic.  A number of intersection points of a last (which is not normalized) skew-section in the marked family with the marked subbundle is even.  A set of intersection points of the section $\psi_2$ with the marked subbundle is represented by an even number of  pairs, points in each pair is covered a corresponding intersection point of the section $\hat{\psi}_2$. Singes in a pair coincides.
The formula (\ref{o})  is proved.

Denote the value of this obstruction  by
 $o(\Psi_2)=4k$, $k \in \Z$.
 Consider the standard immersion  
$f: S^{2^l-1} \looparrowright \R^{2^{l+1}-2}$ with the only self-intersection point. For the immersion $f$ a stable framing is well-defined, this stable skew-framing can be considered as the corresponding stable skew-framing. 
For this skew-framing $\Xi^{st}$ we get:
$o(\Xi^{st})=\pm 2$, where a sign is related of a coorientation of the standard immersion $f$. Take an even $-2k$ algebraic copies of the skew-framed immersion $f$, join this copies 
to the immersion
  $(\phi_2,\Psi_2)$, we get a stably skew-framed immersion  $(\phi_2,\Psi_2)$ with $o(\Psi_2)=0$,
  which is connected in its stably skew-framing cobordism class.
  As the result, we get a skew-framed immersion $\phi_3 : N_3^{2^l-1} \looparrowright \R^{2^{l+1}-2}$
  with an odd number of self-intersection points. Theorem
    \ref{Th1} is proved.

\section{Discussion}
In the section
 \ref{sec3} an element   
 $\eta^{sf}_l \in Imm^{sf}(2^l-1,1)$, $l \ge 3$ is constructed.  
 This element is associated by the Khan-Priddy transfer with the element
 $\eta_l \in \pi_{2^{l+1}+2}(S^{2^l+2}) \simeq \Pi_{2^l-2}$, which is not
 desuspended $5$ times from the stable range into un unstable range.
 
 The same time, for an arbitrary
   $d \ge 0$ there exists $l_0=l_0(d)$ such that for an arbitrary
$l \ge l_0$ an arbitrary element in the stable homotopy group of spheres $\Pi_{2^{l}-2}$ admits a $d$-fold desuspension, as it is proved in  \cite{A1}. In \cite{L-R} it is proved, that if the group $\Pi_{126}$ there exists an element with the Kervaire invariant one then  in the residue class of this element there exists an element, which is not desuspended $11$ times into an unstable range. 

The author is grateful S.A.Melikhov for discussions of the Section \ref{sec3}.

\end{document}